\documentclass{article}

\def\Z{{\rm Z}}
\def\N{{\rm N}}

\def\R{{\rm R}}

\begin{document}

\bigskip

Tsemo Aristide

College Boreal

951, Carlaw Avenue,

M4K 3M2

 Toronto, ON Canada.

tsemo58@yahoo.ca

\bigskip

\centerline{\bf Gerbes for the Chow.}

\bigskip
\bigskip

{\it You can try to stop me, but it wont do a thing no matter what
you do, I'm still gonna be here Through all your lies  and silly
games I'm  still remain the same, I'm unbreakable }

         Michael Joseph Jackson, unbreakable.

\bigskip

\centerline{\bf Abstract.}

Finding coherent relations to define non Abelian cohomology is a
thriller which entertains the mathematical community since fifty
one years. The purpose of this paper is to simplify the attempt to
beat it defined by the author which used the notion of sequences
of fibred categories and to apply the resulting theory to higher
divisors and Chow theory.

\bigskip

\centerline{\bf Introduction.}

\bigskip

The A.B.C of non Abelian cohomology has been created by
Grothendieck and his collaborators in the purpose of giving a
geometric interpretation of characteristic classes. Let $(C,J)$ be
a Grothendieck site, and $L$ a sheaf defined on $(C,J)$, we know
that $H^0(X,L)$ is the set of global sections of $L$, and
$H^1(X,L)$ is the set of isomorphism classes of torsors bounded by
$L$. In his book [5], Giraud  has defined the notion of gerbes
bounded by a sheaf $L$, objects which are classified by
$H^2(X,L)$. Many concrete problems have created the need to
provide a geometric interpretation of higher cohomology classes.
Specialists who wanna be starting something have developed   many
attempts to find a theory which interprets higher cohomology
classes, all of these theories face at this time, the combinatoric
problems which arise when one try to pin out coherence relations
for $n$-gerbes, $n>2$.

In [17],  is developed the notion of sequences of gerbes which
provides a partial answer by attaching to a sequence of fibred
categories endowed with nice properties a cohomology class, this
construction can also be viewed as a geometric interpretation of
the connecting morphism in cohomology. Remark that this approach
gives a complete satisfaction in the geometric study of the Brauer
group as shows [14]. This theory has been successfully applied in
many areas, like in symplectic geometry, where it has enabled to
give new insights on quantization and symplectic fibrations. It
has also been applied to the study of moduli spaces in
differential geometry [19].

The purpose of this paper is to simplify the notion of sequences
of fibred categories studied in [17]. The main tool used here is
the topos of the site of sheaves $Sh(C,J)$ defined on the
Grothendieck site $(C,J)$. In [7], Grothendieck defines on
$Sh(C,J)$ a Grothendieck topology, which can be used to define
notions of varieties and algebraic spaces for any Grothendieck
site (see also [16]). This topology allows us to define here
fibred categories on the basis $(C,J)$ for which the objects of
the fibres are varieties, thus are naturally endowed with a
Grothendieck topology; we study such $2$-sequences of fibred
categories, and apply our results to define and study higher
divisors in algebraic geometry.

\bigskip

 {\bf Notations.}

In this paper all categories are stable by finite limits and
colimits. Let $C$ be a category. We denote by $C/X$, the comma
category of morphisms of $C$ whose target is $X$. Let
$U_{i_1},...,U_{i_p}$ be objects of $C$, we denote by
$U_{i_1...i_p}$, the fiber product of $U_{i_1},...,U_{i_p}$ over
the terminal object of $C$. Let $p:F\rightarrow C$ be a fibred
category (see definition 1, paragraph 2), and a morphism
$h_{i_1...i_p}:e_{i_1...i_p}\rightarrow e'_{i_1...i_p}$ between
objects of the fibre $F_{U_{i_1...i_p}}$, we denote by
$h^{j_1...j_l}_{i_1...i_p}:e^{j_1...j_l}_{i_1...i_p}\rightarrow
{e'}^{j_1...j_l}_{i_1...i_p}$ the restriction of $h_{i_1...i_p}$
between the respective restrictions $e^{j_1...j_l}_{i_1...i_p}$
and ${e'}^{j_1...j_l}_{i_1...i_p}$ of $e_{i_1...i_p}$ and
${e'}_{i_1...i_p}$ to $U_{i_1..i_pj_1..j_l}$.

\bigskip

{\bf 1. Grothendieck topologies, varieties and geometric spaces.}

\medskip

{\bf Definition 1.} Let $C$ be a category, a sieve $S$ defined on
$C$, is a subclass of the class $ob(C)$, of objects of $C$ such
that: if $X\in S$, and  $Y\rightarrow X$ is a morphism of $C$,
then $Y\in S$.

A  Grothendieck topology $J$, defined on the category $C$, is a
correspondence which assigns  to every object $X$ of $C$, a non
empty class of sieves $J(X)$ of  $C/X$ such that:

- Let $S\in J(Y)$, and $f:X\rightarrow Y$ a morphism of $C$, the
pullback $S^f=\{h:Z\rightarrow X, f\circ h\in S\}$ is an element
of $J(X)$.

- A sieve $S$ of $C/X$ is an element of $J(Y)$ if for every
morphism $f:X\rightarrow Y$, $S^f\in J(X)$.

A category equipped with a Grothendieck topology is called a site
or a Grothendieck site. We write it $(C,J)$.

\bigskip

{\bf Definitions 2.} Let $C$ be a category, a presheaf $F$ on $C$
is a functor $F:C^0\rightarrow Set$, where $C^0$ is the opposite
category  of $C$, and $Set$ the category of sets. We denote by
$PreSh(C)$ the category of presheaves defined on the category $C$.

A presheaf on the Grothendieck site $(C,J)$ is called a sheaf if
and only if for every object $X$ of $C$, and every element $S$ of
$J(X)$, $Lim_{Y\rightarrow X\in S}F(Y)=F(X)$.

Let $C$ be a site, a trivial sheaf $F$ defined on $C$, is a sheaf
$F$ such that there exists a set $E$ such that for every object
$X$ of $C$, $F(X)=E$, and the restriction maps are the identity of
$E$.

We denote by $Sh(C,J)$ the category of sheaves on the Grothendieck
site $(C,J)$. We have an full imbedding $C\rightarrow PreSh(C)$
defined by the Yoneda imbedding which associates to the object $X$
of $C$, the presheaf $h_X$ defined by $h_X(Z)=Hom_C(Z,X)$. We say
that the topology is subcanonical if the presheaf $h_X$ is a sheaf
for every $X\in C$. In the sequel, we will consider only
subcanonical topologies. If there is no confusion, we will often
denote $h_X$ by $X$.

\bigskip

{\bf Examples.}

\medskip

Let $Op$ be the category whose objects are open subsets of
${\R}^n, n\in{\N}$, and whose morphisms are local homeomorphisms;
for every object $X$ of $Op$, an element of $J(X)$ is a family of
local homeomorphisms $(h_i:U_i\rightarrow X)_{i\in I}$ such that
$\bigcup_{i\in I}h(U_i)=X$. Every  topological space $T$ defines a
sheaf $h_T$ of $Op$ by assigning to $X$ the set of continuous
maps: $h_T(X)=Hom(X,T).$

\medskip

Let $Aff$ be the category of affine schemes: it is the category
opposite to the category of commutative rings with a unit. We
endow $Aff$ with the etale topology. For every object $X$ of
$Aff$, an element of $J(X)$ is a finite family of etale morphisms
$(h_i:U_i\rightarrow X)_{i\in I}$ such that $\bigcup_{i\in
I}h(U_i)=X$. Every scheme $S$ defines a sheaf $h_S$ of $Aff$ by
assigning to $X$ the set of morphisms of schemes:
$h_S(X)=Hom(X,S)$. See [8], VIII. proposition 5.1.

\medskip

{\bf Definition 3.} Let $(C,J)$ be a site, we say that the
morphism  $F\rightarrow G$ between elements  of $Sh(C,J)$ is a
covering morphism  if and only if  for every object $X$ in $C$,
and every morphism $X\rightarrow G$, the canonical projection
$X\times_GF\rightarrow X$ is a covering sieve of $X$; the family
of morphisms $(F_i\rightarrow G)_{i\in I}$ is a covering family of
$G$ if and only if the morphism $\coprod_{i\in I}F_i\rightarrow G$
is a covering morphism, see [7], p 251-252. These covering
families define on $Sh(C,J)$ a Grothendieck topology (See [7]
proposition 5.4 p. 254).

\medskip

{\bf Definition 4.} Let $C$ be a category, a monomorphism of $C$
is a morphism $f:X\rightarrow Y$, such that for every object $Z$
of $C$, The map $Hom(Z,X)\rightarrow Hom(Z,Y)$ which sends the
element $h\in Hom(Z,X)$ to $f\circ h$ is injective.

Suppose that $C$ is a site, denote by $e$ the final object of
$Sh(C,J)$. The object $X$ of $C$ is an open subset of $e$  if and
only if there exists a monomorphism $i:X\rightarrow e$ (see [1] p.
20 and [7] definition 8.3 p. 421). The morphism $i$ is called an
open immersion.

 The object $U$ of $C/X$ is an open subset of $X$, if and only if it is an open subset
 of the final object of $C/X$ for the induced topology.

\medskip

{\bf Definition 5.} Let $(C,J)$ be a site, we suppose that for
every object $X$ of $C$,  every open subset $f:U\rightarrow X$, of
$C/X$ is  contained in a sieve of $X$, a geometric space is a
sheaf $F$ of $(C,J)$ such that:

 There exists a family $(U_i)_{i\in I}$ of objects of $C$ and a
sieve $p:\coprod_{i\in I}{U_i}\rightarrow F$ of $F$, for the
Grothendieck topology on $Sh(C,J)$.

The family $(U_i)_{i\in I}$ is called an atlas.

Let $p_i:{U_i}\rightarrow F$ be the composition of the canonical
imbedding ${U_i}\rightarrow \coprod_{i\in I}{U_i}$ and $p$. If for
every $i$, the map $p_i$ is an open immersion, then $F$ is called
a variety.

\bigskip

{\bf Examples.}

A geometric space in $Op$ is defined by a sheaf on $Op$, and a
covering morphism $p:\coprod_{i\in I}U_i\rightarrow F$. In
particular a topological manifold is a geometric space, it is in
fact a variety.

\medskip

A geometric space $F$ in $Aff$ is defined by a sheaf on $Aff$, and
a covering morphism $p:\coprod_{i\in I}U_i=Spec(A_i)\rightarrow
F$. In particular a scheme is a geometric space, in fact it is a
variety.

\bigskip

Let $F$, be a geometric space, suppose that the covering
$p:\coprod_{i\in I}U_i\rightarrow F$ is $1$-connected; (this is
equivalent to saying that for every $i\in I$, every sheaf on $U_i$
is trivial). The pullback $F_i$ of $F$ by $p_i$ is trivial. Let
$F_i^j$ be the pullback of $F_i$ on $U_i\times_FU_j$ by the
projection $U_i\times_FU_j\rightarrow U_i$. There exists an
isomorphism $g_{ij}:F_j^i\rightarrow F_i^j$.  The morphism
$c_{ijk}=g_{ki}^jg_{ij}^kg_{jk}^i$ is an automorphism of
$F_k^{ij}$, the restriction of $F_k$ to
$U_i\times_FU_j\times_FU_k$. Which can be identified with an
automorphism of $F(U_i\times_FU_j\times_FU_k)$.

The family $c_{ijk}$ verifies the relation:
$c^j_{ikl}g^{ij}_{lk}c^l_{ijk}g^{ij}_{kl}=c^k_{ijl}c^i_{jkl}$, we
say that it is a (non commutative) $2$-cocycle.

\medskip

Suppose that $F$ is a variety; recall that this is equivalent to
saying that ${U_i}\rightarrow F$ is an open immersion for every
$i$. Thus $F_{i}$ is  the restriction of $F$ to $U_i$. Let
$h_i:F_i\rightarrow {E_i}$, the trivialization of the restriction
of $F$ on $U_i$, on $U_i\times_FU_j$, we can define the map
${h_i^j}^{-1}\circ{h^i_j}$ which is an automorphism of $E_{ij}$,
the fibre of the restriction of $F$ to $U_i\times_FU_j$. We have
the relation: $c_{ik}^j=c_{ij}^kc_{jk}^i$

\medskip

{\bf Definition 6.} Let $C$ be a category, $X,P$ objects of $C$, a
$P$-point of $X$ is a morphism $x:P\rightarrow X$.

In a category stable by finite limits, and colimits, a group
object (See [1] p.35) is defined by:

 - An object $G$ endowed with a morphism
$p:G\times G\rightarrow G$ called the product, which is
associative,

- The neutral element, which is a global point, that is a morphism
$e:1\rightarrow G$, (where $1$ is the final object).

- The inverse is a morphism $i:G\rightarrow G$.

This data must satisfy the following conditions:

 Let $x,y:P\rightarrow G$ be two $P$-points of $G$, by the universal
property of the product, $x$ and $y$ define  a morphism
$(x,y):P\rightarrow G\times G$, we write $p\circ (x,y)=xy$, we
must have $x(yz)=(xy)z$.

 Let $i(x)$ be the point $i\circ x$, the fact
that $i$ is the inverse map is equivalent to saying that
$p(x,i(x))$ is the composition of the unique map $p_P:P\rightarrow
1$ with $e$, we must also have $p(x,e\circ p_P)=(e\circ p_P,x)=x$.

An action of the group object $G$ on $X$ is defined by a morphism
$A:G\times X\rightarrow X$. The universal property of the product
and the action induces a morphism: $h_{T,X}:Hom(T,G)\times
Hom(T,X)\rightarrow Hom(T,X)$; for every points $g,g':T\rightarrow
G$, and $x:T\rightarrow X$, $(gg')x=g(g'x)$.

The action is free if and only if for every $T$-point, $g$ of $G$,
$h_{T,X}(g,.)$ is injective.

 Remark that by using the Yoneda imbedding,  if $h_G$ is a
  group object of $Sh(C,J)$, then $G$ is also a group object of $C$.

\bigskip

{\bf Proposition 1.} {\it Let $G$ be a group object of $Sh(C,J)$
which acts freely on the geometric space $X$, and such that for
every object $Y\in C$, the projection $Y\times G\rightarrow Y\in
J(Y)$ then the sheaf $h_X/G$ is a geometric space.}

\medskip

{\bf Proof.} Let $({U_i})_{i\in I}$ be an atlas of $X$; Denote by
$p'_i$ the composition of $p_i$ with the map $X\rightarrow X/G$
(see definition 5 for $p_i$). Then $(U_i,p'_i)_{i\in I}$ defines
an atlas of $X/G$. To show this, firstly we consider an object $Z$
of $C$, and a morphism $h:h_Z\rightarrow X/G$. The pullback of
$X\rightarrow X/G$ by $h$ is $h_Z\times h_G$; to show this,
consider an element of $Hom(T,Z)\times_{Hom(T,X)/G}Hom(T,X)$ which
is defined by an element $u$ of $Hom(T,Z)$, and an element $v$ of
$Hom(T,X)$ which have the same image in $Hom(T,X)/G$. The elements
of $Hom(T,Z)\times_{Hom(T,X)/G}Hom(T,X)$ whose image by the first
projection is $u$ are of the form $(u,gv), g\in Hom(T,G)$. Since
the action of $G$ is free, we deduce that
$h_Z\times_{X/G}h_X=h_Z\times h_G$.

Since $(U_i)_{i\in I}$ is an atlas of $X$, the pullback of
$h_{Z\times G}\rightarrow X$ by $\coprod_{i\in I}
h_{U_i}\rightarrow X$ is in $J(Z\times G)$, since the map $Z\times
G\rightarrow Z\in J(Z)$. We deduce that $(U_i,p'_i)_{i\in I}$ is
an atlas of $X/G$.

\bigskip

{\bf 2. Sheaves of categories $(2,2)$-gerbes, $(2,1)$-gerbes,
$(1,2)$-gerbes.}

\medskip

Let $(C,J)$ be a category equipped with a Grothendieck topology,
and $p:F\rightarrow C$ a functor. For every object $X$ of $C$, we
denote $F_X$, the subcategory of objects of $F$, such that for
every object $x$ of $F_X$, $p(x)=X$. A morphism $h:x\rightarrow
x'$ between objects of $F_X$ is an element of $h\in Hom_F(x,x')$,
such that $p(h)=Id_X$. The category $F_X$ is called the fibre of
$X$. Let $f:X\rightarrow Y$ be a morphism of $C$, and $x\in F_X,
y\in F_Y$. We denote by $Hom_f(x,y)$ the subset of the set of
$Hom_F(x,y)$ such that for every element $h\in Hom_f(x,y)$,
$p(h)=f$.

\medskip

{\bf Definitions 1.} The morphism $h\in Hom_f(x,y)$ is Cartesian
if and only if for every element $z\in F_X$, the canonical map
$Hom_{Id_X}(z,x)\rightarrow Hom_f(z,y)$ which sends $l\rightarrow
h\circ l$ is bijective.

We say that the category $F$ is a fibred category over $C$, if and
only if for every morphism $f:X\rightarrow Y$ in $C$, and every
element $y\in F_Y$, there exists a Cartesian morphism
$c_f:x\rightarrow y$ such that $p(c_f)=f$.

Examples: the forgetful functor $C/X\rightarrow C$, which sends
$Y\rightarrow X$ to $Y$, is Cartesian, as well as its restriction
to any sieve of $X$. Let $p:F\rightarrow C$ and $p':F'\rightarrow
C$ be Cartesian functors, we denote by $Cart(F,F')$ the class of
morphisms between $F$ and $F'$ such that for
 every element $h\in Cart(F,F')$, we have $p'\circ h=p$, and $h$
 sends Cartesian morphisms to Cartesian morphisms.

 \medskip

{\bf Definition 2.} Let $p:F\rightarrow C$ be a Cartesian functor.
We say that $F$ is a sheaf of categories, if and only if:

- For every sieve $R\in J(X)$, the forgetful functor
$Cart(E/X,F)\rightarrow Cart(R,F)$ is an equivalence of
categories.

We say that the sheaf of categories is connected if for every
object $X$ of $C$, there exists a sieve $R\in J(X)$, such that for
every morphism $Y\rightarrow X\in R$, $F_Y$ is not empty, and the
objects of $F_Y$ are isomorphic each others. We are going to study
only connected sheaves of categories here.

Let $f:X\rightarrow Y$ be a morphism of $C$, and $y$ an object of
$F_Y$, a restriction map of $f$  is a Cartesian map
$c_f:x\rightarrow y$, we say often that $x$ is a restriction of
$y$.

\medskip

Suppose that the topology of $C$ is generated by the family
$(U_i)_{i\in I}$, we can assume that for every $i\in I$, the
object of the fibre of $U_i$ are isomorphic each other. Choose an
object $x_i\in F_{U_i}$, on $U_{ij}$, there exists a morphism
$g_{ij}:x^i_j\rightarrow x^j_i$; the morphism
$c_{ijk}=g_{ki}^jg_{ij}^kg_{jk}^i$ is an automorphism of
$x_k^{ij}$. Which satisfies the relation:

$$
c^j_{ikl}g^{ij}_{lk}c^l_{ijk}g^{ij}_{kl}=c^k_{ijl}c^i_{jkl}
$$

We have seen that geometric spaces satisfy the condition above, in
fact, there are examples of sheaves of categories.

\medskip

{\bf Definition 3.} A sheaf of categories on the Grothendieck site
$(C,J)$ is a gerbe, if and only if:

 there exists a sheaf $L$ on $(C,J)$ such that for every object $x\in F_X$,
$Aut_{Id_X}(x)\simeq L(U)$, and this identification commutes with
morphisms between objects and with restrictions. We say that the
sheaf of categories $p:F\rightarrow C$ is bounded by $L$.

\bigskip

{\bf Definition 4.} Suppose that the site $(C,J)$ has a final
object $e$; a gerbe is trivial if and only if it has a global
section. This is equivalent to saying that the fibre $F_e$ is not
empty.

A global section is called a torsor; equivalently a torsor is a
gerbe $p:F\rightarrow C$ such that for every object $X$ of $C$,
the fibre $F_X$ contains a unique object.

\bigskip

 {\bf Definition 5.}
Let $(C,J)$ be a Grothendieck site. Consider a variety $X$ defined
on $(C,J)$. An $n$-sequence of fibred categories over $X$, is a
sequence of functors $p_n:F_n\rightarrow F_{n-1}...
p_1:F_1\rightarrow C/X=F_0$ which satisfies the following
conditions:

- The functors $p_l, l=1,..n$ are fibred categories.

- For every object $U$ of $F_l$, the fibre ${F_{l+1}}_U$ is a
category whose objects are varieties of $C$, and its morphisms are
morphisms of varieties over $U$.

\bigskip

To define the notion of $n$-sequence of gerbes, we are going to
associate firstly, to an automorphism above the identity of a
gerbe bounded by a commutative sheaf, a $1$-cocycle. Let $h$ be an
automorphism above the identity of the gerbe $p:F\rightarrow C/X$.
Let $(U_i)_{i\in I}$ be $1$-connected cover of $X$ (This is
equivalent to saying that the restriction of every sheaf to $U_i$
is trivial), let $x_i$ be an object of $F_{U_i}$, there exists an
arrow $l_i:x_i\rightarrow h(x_i)$. Let $u_{ij}:x_j^i\rightarrow
x_i^j$ be a connecting morphism, on $U_{ij}$ we have the morphism
$h_{ji}={l_j^i}^{-1}\circ h(u_{ij})^{-1}\circ l_i^j\circ u_{ij}$
of $x_j^i$. We are going to show that
$h^i_{kj}u_{kj}h^k_{ji}u_{jk}=h^j_{ki}$:

$h^i_{kj}u_{kj}h^k_{ji}u_{jk}={l_k^{ij}}^{-1}
h(u_{kj}^i)^{-1}{l_j^{ik}}u^i_{jk}u^i_{kj}{l_j^{ik}}^{-1}h(u^k_{ij})^{-1}
{l_i^{jk}}u^k_{ij}u^i_{jk}={l_k^{ij}}^{-1}h(u_{ij}^iu^k_{jk})^{-1}l_i^{jk}u^k_{ij}u^i_{jk}$.
By writing that $c_{ijk}=u_{ki}^ju_{ij}^ku_{jk}^i$, we obtain
that:

$h^i_{kj}u_{kj}h^k_{ji}u_{jk}={l_k^{ij}}^{-1}h^{-1}(u^j_{ik}c_{ijk})l_i^{jk}u_{ik}c_{ijk}$.
Since the group $L$ is commutative and $h$ commutes with morphisms
between objects,
$h^{-1}(u^j_{ik}c_{ijk})l_i^{jk}u_{ik}c_{ijk}=c_{ijk}^{-1}h^{-1}(u^j_{ik})l_i^{jk}c_{ijk}u_{ik}
=h^{-1}(u^j_{ik})l_i^{jk}u_{ik}$. This implies that
$h^i_{kj}u_{kj}h^k_{ji}u_{jk}=h^j_{ki}$.

The cohomology class of the cocycle that we have just defined
doesn't depend of the choices made. Suppose that we fix the $x_i$,
but replace $l_i$ by $l'_i$, then there exists $u_i\in L(U_i)$
such that $l'_i=u_il_i$, and $h_{ij}$ is replaced by
${u^i_j}^{-1}h_{ij}u_i^j$.

Suppose that we replace $x_i$ by $x'_i$, let $v_i:x'_i\rightarrow
x_i$ be a connecting morphism, $h(v_i)^{-1}l_iv_i$ is a connecting
morphism $l'_i$ between $x'_i$ and $h(x'_i)$,
$u'_{ij}={v_i^j}^{-1}u_{ij}v_j^i$ is a connecting morphism between
${x'_j}^i$ and ${x'_i}^j$. We can write
${{l'}^i_j}^{-1}h^{-1}(u'_{ij}){l'_i}^ju'_{ij}=$

$({h(v_j^i)^{-1}l_j^iv_j^i})^{-1}h({v_i^j}^{-1}u_{ij}v_j^i)^{-1}
h(v^j_i)^{-1}l_i^jv_i^j{v_i^j}^{-1}u_{ij}v_j^i={v_j^i}^{-1}h_{ij}v_j^i=h_{ij}$
since the elements of the band commute with morphisms between
objects.

\medskip

{\bf Definition 6.} A $n$-sequence of fibred categories
$p_n:F_n\rightarrow F_{n-1} ...p_1:F_1\rightarrow C/X=F_0$ is a
$n$-sequence of gerbes, if and only if:

- For every object $U$ of $F_{n-2}$, and $e_U$ of ${F_{n-1}}_U$,
the fiber ${F_n}_{e_U}$ is a gerbe bounded by a sheaf $L_{e_U}$
defined  on $C/e_U$.

- Let $U$ be an object of $F_l$. There exists a cover $(U_i)_{i\in
I}$ of $U$, such that for every object $e_i,e'_i$ of
${F_{l+1}}_{U_i}$, and there exists an isomorphism between
${F_{l+2}}_{e_i}$ and ${F_{l+2}}_{e'_i}$.

- There exists  a commutative sheaf $L$ on $C$  called the band
such that the trivial automorphisms (those corresponding to
trivial bundles) of ${F_n}_{e_U}$ are the sections of $L$.

\bigskip

{\bf The classifying $4$-cocycle.}

\medskip

In the sequel, we will consider only $2$-sequences of fibred
categories that we call also $(2,2)$-gerbes, the general situation
will be studied in a forthcoming paper.

We are going to associate a $4$-cocycle to a $2$-sequence
$p_2:F_2\rightarrow p_1:F_1\rightarrow X$ bounded by a sheaf of
commutative groups $L$.

Let $(U_i)_{i\in I}$ be a cover of $X$, and $x_i$ an object of
${F_1}_{U_i}$, we denote by $g_{ij}:x_j^i\rightarrow x_i^j$ a
connecting morphism. The morphism
$c_{ijk}=g_{ki}^jg_{ij}^kg_{jk}^i$ is an automorphism of
$x_k^{ij}$. Let $U$ be an object of $C/x_k^{ij}$, for every object
$U'\in {F_2}_U$ we can lift the pullback of $c_{ijk}$ by
$U\rightarrow x_{ij}^k$, to a Cartesian morphism $U"\rightarrow
U'$. If $h':U'\rightarrow V'$ is a morphism above the morphism
$h:U\rightarrow V\rightarrow x_k^{ij}$ between objects of
$C/x_k^{ij}$, we can lift the pullback of $h$ by $c_{ijk}$ to a
morphism $h":U"\rightarrow V"$ in such a way that
$h":U"\rightarrow V"\rightarrow V'$ coincide with $U"\rightarrow
h':U'\rightarrow V'$. This shows that the correspondence which
associates $U"$ to $U'$ defines  an automorphism  $c'_{ijk}$, of
the gerbe ${F_2}_{x^{ij}_k}$. (See also Giraud [6] Scholie 1.6
p.3)

On $U_{ijkl}$, we have the morphisms $c^i_{jkl}, c^j_{ikl},
c^k_{ijl}, u_{lk}^{ij}c^l_{ijk}u_{kl}^{ij}={c'}_{ijk}$ of
$x_l^{ijk}$. The automorphism
${{c'}^l_{ijk}}^{-1}{c^j_{ikl}}^{-1}c_{ijl}^kc_{jkl}^i$ is an
automorphism above the  identity of $x_l^{ijk}$. We identify it
with an element $c_{ijkl}\in C^1(x_l^{ijk},L)$ up to a boundary.
The cohomology class of the Cech boundary $c_{ijklm}$ of
$c_{ijkl}$ is trivial. Thus we can identify $c_{ijklm}$  with an
element of $L(U_{ijklm})$ The family $c_{ijklm}$ is the
classifying $4$-cocycle of the $(2,2)$-gerbe.

\bigskip

{\bf $(2,1)$-gerbe, and $(1,2)$-gerbe.}

\medskip

We will often need a particular $2$-sequence of gerbes:

\medskip

{\bf Definition 7.} A gerbe-torsor or a $(2,1)$-gerbe is a
$(2,2)$-gerbe $p_2:F_2\rightarrow p_1:F_1\rightarrow X$, such
that:

- For every object $U$ of $C/X$, there exists a covering
$(U_i)_{i\in I}$ of $X$ such that for every object $e_{U_i}$ of
${F_1}_{U_i}$, the category ${F_2}_{e_{U_i}}$ is a trivial gerbe
over $e_{U_i}$.

- There exists a sheaf $L$ such for every global section  $V$ of
${F_2}_{e_{U_i}}$, we can identify $Aut_{e_{U_i}}(V)$, the group
of automorphisms of $V$ over the identity of $e_{U_i}$ with
$L(U_i)$, and this identification commutes with morphisms between
objects and with restrictions.

\medskip

{\bf The classifying cocycle of a $(2,1)$-gerbe.}

\medskip

We are going to associate to a gerbe-torsor,  a $3$-cocycle
defined as follows:

 Let $x_i$ be an object of ${F_1}_{U_i}$, and
 $u_{ij}:x_j^i\rightarrow x_i^j$ a morphism, we can define the
 cocyccle $c_{ijk}=u_{ki}^ju_{ij}^ku_{jk}^i$ of $x_k^{ij}$. Since
 the gerbe ${F_2}_{x_k^{ij}}$ is trivial, we can pick $V$, a global
 section  over $x_k^{ij}$, in this situation, let  $c'_{ijk}$  be
 a morphism of $V$ above $c_{ijk}$.  The Cech boundary of $c'_{ijk}$
 is an automorphism  above the identity of $x_l^{ijk}$ that we identify with an
 element of $L(U_{ijkl})$. The family of morphisms $c_{ijkl}$
 defines a $3$-cocycle which is $L$-valued.

 \medskip

 {\bf Remark.}

 Suppose that the cohomology class of $c_{ijkl}$ is zero, thus up
 to a boundary, we can assume that $c_{ijkl}=0$. This is
 equivalent to saying that $c'_{ijk}$ is the classifying cocycle
 of a gerbe $p':F'\rightarrow C/X$, such that for every object $U$
of $C/X$, the objects of ${F'}_U$ are torsors over $e_U$, where
$e_U$ is an object of ${F_1}_U$. The classifying cocycle of the
$(2,1)$-gerbe can be viewed as an obstruction to obtain such a
gerbe, that this to reduce the trivial gerbe ${F_2}_{e_U}$ to a
torsor.

\medskip

Let $p_2:F_2\rightarrow p_1:F_1\rightarrow C/X$ be a
$(2,2)$-gerbe, we are going to associate to it a $(2,1)$-gerbe
$p'_2:F'_2\rightarrow p_1:F_1\rightarrow C/X$ defined as follows:
Let $U$, be an object of $C/X$, and $e_U$ an object of ${F_1}_U$,
the trivial gerbe ${F'_2}_{e_U}$ is the gerbe whose objects are
the automorphisms of ${F_2}_{e_U}$ above morphisms of $e_U$.
Remark that the classifying cocycle of this $(2,1)$-gerbe is the
class $c_{ijkl}$ that we have used to define the classifying
cocycle of $(p_2,p_1)$. Thus $(p_2,p_1)$ is trivial if and only if
$(p'_2,p_1)$ is trivial. This allows to interpret a $(2,2)$-gerbe
as a geometric obstruction.

\bigskip

{\bf Definition 8.} A $(1,2)$-gerbe is a $(2,2)$-gerbe
$p_2:F_2\rightarrow p_1:F_1\rightarrow C$, such that the gerbe
$p_1:F_1\rightarrow C$ is trivial.

\bigskip

{\bf The classifying $3$-cocycle of a $(1,2)$-gerbe.}

\medskip

Suppose that the covering $(U_i)_{i\in I}$ is a good covering, and
 $x_i={F_1}_{U_i}$ is isomorphic to the trivial $L(U_i)$-torsor
$t_i$ on $U_i$. Let $h_i:x_i\rightarrow t_i$ be an isomorphism,
consider the automorphism $c_{ij}={h_i^j}^{-1}h_j^i$ of $x_{ij}$.
Let $U$ be an object of $C/x_{ij}$. We can lift $c^i_{jk}$ to an
automorphism ${c'}^i_{jk}$ of ${F_2}_{x_{ijk}}$. The Cech boundary
of ${c'}^i_{jk}$ is an automorphism of the gerbe ${F_2}_{x_{ijk}}$
above the identity that we identify with an element of $c_{ijk}\in
C^1(U_{ijk},L)$ up to a boundary. The cohomology class of the
boundary $c_{ijkl}$ of $c_{ijk}$ is trivial. We can thus identify
$c_{ijkl}$ to an element of $L(U_{ijkl})$. The family $c_{ijkl}$
is the classifying cocycle of the $(1,2)$-gerbe.

\bigskip

{\bf Examples: The  lifting obstruction.}

\bigskip

 Let $(C,J)$ be a site, and
$0\rightarrow L\rightarrow M\rightarrow N\rightarrow 0$ an exact
sequence of commutative sheaves defined on $C$. It defines the
following exact sequence in cohomology:

$$
H^n(X,L)\rightarrow H^n(X,M)\rightarrow H^n(X,N)\rightarrow
H^{n+1}(X,L)
$$

 Let $[c^n]$ be an element of $H^n(X,N)$, represented by the
 $n$-cocycle $c^n$ of the sheaf $N$. A natural problem is to find obstructions to lift
$[c^n]$ to a $H^n(X,M)$.

Recall the construction of the boundary operator
$H^n(X,N)\rightarrow H^{n+1}(X,L)$. Let $(U_i)_{i\in I}$ be a good
cover of $(C,J)$, the restriction of $M$ and $N$ on $U_i$ are
trivial. This implies the existence of a global section $b^n_i\in
M(U_{i_1..i_{n+1}})$ over $c^n_i$, the restriction of $c_n$ to
$U_{i_1..i_{n+1}}$. We can write the boundary $c^{n+1}$ of the
chain $b^n_i$ it is an $L$-cocycle whose cohomology class is the
image of $[c^n]$ by the boundary operator.

The cohomology classes of the $n+1$-$L$-cocycles which are in the
image of the connecting morphism $H^n(X,N)\rightarrow
H^{n+1}(X,L)$ are in bijection with the quotient of $H^n(X,N)$ by
the image of the morphism $H^n(X,M)\rightarrow H^n(X,N)$.

\bigskip

If $n=0$, $H^0(X,N)$ classifies the global sections of the sheaf
$N$, and $H^1(X,L)$ the $L$-torsors, we obtain that isomorphism
classes $L$-torsors whose classifying cocycles are in the image of
the connecting morphism $H^0(X,N)\rightarrow H^1(X,L)$ are in
bijection with the quotient of $H^0(X,N)$ by the image of the
morphism $H^0(X,M)\rightarrow H^0(X,N)$.

\medskip

If $n=1$, $H^1(X,N)$ classifies the torsors of the sheaf $N$, and
$H^2(X,L)$ the $L$-gerbes, we obtain that isomorphism classes of
$L$-gerbes whose classifying cocycles are in the image of the
connecting morphism $H^1(X,N)\rightarrow H^2(X,L)$ are in
bijection with the quotient of $H^1(X,N)$ by the image of the
morphism $H^1(X,M)\rightarrow H^1(X,N)$.

\medskip

Let $p^i_2:F^i_2\rightarrow p^i_1:F^i_1\rightarrow C/X, i=1,2$ be
two $(2,1)$-gerbes such that $p^i_1$ is a gerbe bounded by $N$,
and $(p^i_1,p^i_2)$ by $L$. We say that they are isomorphic if and
only if the respective $3$-cohomology classes associated to these
gerbes are equal.

\medskip

{\bf Proposition 1.} {\it Suppose that the morphism of sheaves
$M\rightarrow N$ has local sections. Then the isomorphism classes
of $(2,1)$-gerbes $p_2:F_2\rightarrow p_1:F_1\rightarrow C/X$ such
that $p_1$ is a gerbe bounded by $L$, and $(p_1,p_2)$ by $N$,
whose classifying cocycles are in the image of the connecting
morphism $l_2:H^2(X,N)\rightarrow H^3(X,L)$ are in bijection with
the quotient of $H^2(X,N)$ by the image of the morphism
$H^2(X,M)\rightarrow H^2(X,N)$.}

\medskip

{\bf Proof.} We need only to construct for every cohomology class
$[c_3]\in H^3(X,L)$ in the image of the connecting morphism
$l_2:H^2(X,N)\rightarrow H^3(X,L)$, a $(2,1)$-gerbe classified by
$[c_3]$. Set $[c_3]=l_2([c_2])$. Let $p:F\rightarrow C$ be an
$N$-gerbe bounded by $[c_2]$. We can suppose that for every object
$U$ of $C$, the objects of the fiber $F_U$ are $N$-torsors. Let
$e_U$ be an object of $F_U$, we define ${F_2}_{e_U}$  to be the
category whose objects are $M$-torsors $p_U:V_U\rightarrow e_U$
whose quotient by $L$ is $e_U$. A morphism between two objects of
${F_2}_{e_U}$ is a morphism of $M$-bundles which projects to the
identity of $e_U$.

 If $(U_i)_{i\in I}$ is a good cover of
$C$, and $e_i$ an object of $F_{U_i}$.  The objects of
${F_2}_{e_i}$ are isomorphic; they are trivial bundles since the
map $M\rightarrow N$ has local sections.

The projection $F_2\rightarrow F_1$ is the projection which sends
the $M$-bundle $V_U\rightarrow e_U$ to $e_U$, this projection is
Cartesian. Let $p:V'_{U'}\rightarrow e_{U'}$ be an element of
${F_2}_{e_{U'}}$, an $f:e_U\rightarrow e_{U'}$ a morphism, the
pullback of $p$ by $f$ is a Cartesian morphism above $f$.

\medskip

{\bf Remark.}

Let $p:F\rightarrow C/X$ and $p':F'\rightarrow C/X$ two gerbes
bounded by $N$, we can define the summand $F+F'$ of $F$ and $F'$:
The objects of $(F+F')_U$ are sum of $N$-bundles $e_U$ and $e'_U$,
where  $U$ is an object of $C/X$, and $e_U$ (resp. $e'_U$) an
object an object of $F_U$ (resp. $F'_U$).

Consider the $(2,1)$-gerbes $p_2:F_2\rightarrow p_1:F_1\rightarrow
C/X$ and $p'_2:{F'_2}\rightarrow p'_1:{F'_1}\rightarrow C/X$ whose
classifying cocycle are image of $l_2$; they are isomorphic if and
only if there exists a $N$-gerbe $F_1"$ whose classifying cocycle
is in the image of $H^2(X,M)\rightarrow H^2(X,N)$ such that
$F_1=F'_1+F"_1$.

\bigskip

{\bf 3. Applications to algebraic geometry: Chow groups and higher
divisors.}

\medskip

In the sequel, $X$ will be  a quasi-projective variety of
dimension $n$ defined on the field $k$, $L_X$ the sheaf of non
zero rational functions  defined on $X$. We endow $X$ with the
Zariski topology. Let $U$ be an open subset of $X$, and $f\in
L_X(U)$, we denote by $(f)$ the principal divisor associated to
$f$. The multiplicative group $L_X(U)$ is a ${\Z}$-group, for the
action defined by $(a,f)\rightarrow f^a, a\in{\Z}, f\in L_X(U)$.
Let $h$ be an element of $L_X^l(U)$, $h=(h_1,...,h_l)$, where
$h_i={{a_i}\over{b_i}}, i=1,...,l$ and $a_i, b_i$ are regular
functions. Denote by $CH^l_X(U)$ the linear subspace generated by
the set of irreducible closed subvarieties of $U$ of codimension
$l$ which are local complete intersections; We define
$ch_l(U):L_X^l(U)\rightarrow CH_X^l(U)$ which sends $h$
 to the intersection product $(a_1-b_1)...(a_l-b_l)\in CH^l_X(U)$.
 Remark that the theorem 1 V.21 of Serre [15] describes  the
 elements of the image of $ch_l(U)$ as complete intersections codimension $l$ subvarieties,
 since it implies that if a component of $(a_i-b_i)$ and a component
 of $(a_j-b_j)$ do not intersect properly, their coefficient in
 $(a_i-b_i).(a_j-b_j)$ is zero.

 The  map $ch_l(U)$ is $l$-multilinear for the multiplicative structure,
 it thus factors by a linear map $ch'_l(U):L_X(U)^{\otimes l}\rightarrow CH_X^l(U)$ which
 factors by the quotient map $L^{\otimes l}_X(U)\rightarrow M^l_X(U)$, where $M^l_X(U)$
is the symmetric functions in $l$-variables on $L_X(U)$ for the
multiplicative structure, that is is the quotient of
$L^{l\otimes}_X(U)$ by its subset generated by elements
$(x_1\otimes..\otimes x_l)-\sigma(x_1\otimes..\otimes x_l),
\sigma\in S_l$.

   Since the element of $CH_X(U)$ are local complete
   intersections, for each integer $l$,
    we have an exact sequence of sheaves:

$$
1\rightarrow Z_X(l)\rightarrow M_X^l\rightarrow CH_X^l\rightarrow
1.\leqno{(1)}
$$

Where $Z_X(l)$ is the kernel of the morphism $M^{l}_X\rightarrow
CH^l_X$; we deduce the existence of the following exact sequence
in cohomology:

$$
H^p(X,Z_X(l))\rightarrow H^p(X,M_X^l)\rightarrow
H^p(X,CH_X^l)\rightarrow H^{p+1}(X,Z_X^X(l)).\leqno{(2)}
$$

 Let $p=0,1,2,3$, we define a $p$-gerbe bounded by the sheaf $L$,
 to be a global section of $L$ if $p=0$, a
$L$-bundle if $p=1$, a $L$-gerbe if $p=2$, and a $(2,1)$-$L$-gerbe
if $p=3$. In the sequel $p$ is an integer equal to $0,1$ or $2$.
This restriction is
 due to the fact that for $n>3$, we cannot provide at this time a geometric
 interpretation of this notion.

\medskip

{\bf Definition 1.}
 A  $(p,l)$-Cartier divisor, is defined by a
$p$-chain $(U_{i_1..i_{p+1}},f_{i_1..i_{p+1}})$ of sections of
$M_X^{l}$  such that the image of the Cech boundary
$d(f_{i_1..i_{p+1}})\in Z_X(l)(U_{i_1..i_{p+2}})$. This boundary
is thus the classifying cocycle of a $p+1$-$Z_X(l)$-gerbe
$A(p,l)$. We can also define $p$-gerbe $B(p,l)$ bounded by
$M_X^l/Z_X(l)$ whose classifying cocycle is defined by the classes
of $f_{i_1..i_{p+1}}$ in
$M_X^l(U_{i_1..i_{p+1}})/Z_X(l)(U_{i_1..i_{p+1}})$.

 \medskip

{\bf Remarks.}

 The classifying cocycle of the $p+1$-gerbe $A(p,l)$,
  is the image of the classifying cocycle of $B(p,l)$ by the connecting morphism
  $H^p(X,M_X^l/Z_X(l))\rightarrow
H^{p+1}(X,Z_X(l))$. By comparing $(2)$ with the exact sequence
$H^p(X,Z_X(l))\rightarrow H^p(X,M_X^l)\rightarrow
H^p(X,M_X^l/Z_X(l))\rightarrow H^{p+1}(X,Z_X(l))$ deduced from the
exact sequence $1\rightarrow Z_X(l)\rightarrow M_X^l\rightarrow
M_X^l/Z_X^l$. We deduce that $H^p(X,M_X^l/Z_X(l))$ is isomorphic
to $H^p(X,CH_X^l)$. The isomorphism classes of the $(p,l)$-Cartier
divisors is the quotient $H^p(X,CH_X^l)$ by the image of the
morphism $h_{p,l}: H^p(X,M_X^l)\rightarrow H^p(X,CH_X^l)$. The
elements of $H^p(X,CH_X^l)$ are called $(p,l)$-Weil divisors. Two
$(p,l)$-Weil divisors $D_W$ and $D'_W$ are equivalent if and only
if $D_W-D'_W$ is an element of the image of $h_{p,l}$.

The $p+1$-chain $d(f_{i_1..i_{p+1}})$ is a boundary of elements of
$M_X^l$. Thus correspond to a trivial $M_X^l$ bundle if $p=0$, a
trivial $M_X^l$-gerbe if $p=1$, and a trivial $M_X^l$-$2$-gerbe if
$p=2$, (See Brylinski-McLaughin  [3] for the definition of
$2$-gerbe).

 If $p=0$, and $l=1$ a $(0,1)$-Cartier (resp. a $(0,1)$-Weil divisor) divisor is nothing but a
Cartier divisor (resp. a Weil divisor) in the classical sense. Two
$(0,1)$-Weil divisors are equivalent if and only if they are
equivalent in the classical sense.

 More generally two $(0,l)$-divisors which are equivalent
are rationally equivalent: this follows from the following
argument: let $D_W$ and $D'_W$ be two Weil divisors, suppose that:
$D'_W=D_W+(a_1)...(a_l)$, where $a_i$,...,$a_l$ are regular
functions. Then $D'_W-D_W$ is a principal divisor of
$(a_1)...(a_{l-1})$, it follows from Hartshorne [9] p. 426, that
$D_W$ and $D'_W$ are rationally equivalent.

  Suppose  that $X$ is an affine variety $X$; since the sheaf of rational functions $L_X$ is constant, we
deduce that $H^p(X,M^l_X)=0,p>0$, and $H^0(X,M_X)=K(X)$ the field
of rational functions of $X$. This implies that
$H^p(X,CH^l_X)=H^{p+1}(X,Z_X(l))$.

Suppose that $l=1$, then $Z_X(1)=O_X^*$ the sheaf of invertible
regular functions, we have: $H^0(X,O_X^*)=O_X^*(X)$,
$H^1(X,O_X^*)=Pic(X)$ the Picard group of $X$, and
$H^p(X,O_X^*)=0$ if $p>1$.

\bigskip

{\bf The Cartier divisor associated to a local complete
intersection subvariety.}

\medskip

 Let $Y$ be a closed subvariety of the quasi-projective variety $X$ of codimension $l$ which is a local
 complete intersection, consider an open cover $(U_i)_{i\in I}$ by affine subsets,
 such that $U_i\cap Y$ is the locus $l$ functions
 $(f_i^1,...,f_i^l)$ which defines the element $F_i\in M_X^l(U)$
 obtained by projecting the image of $(f_i^1,...,f_i^l)$. The
 element $h_{ij}=F_j-F_i$ is in $Z_X^l(U_i\times_XU_j)$.

 \medskip

 {\bf Proposition 1.}
 {\it The element $h_{ij}$ is a $1$-Cech $Z_X^l$ cocycle. If $Y$ is irreducible, its
 cohomological class vanishes if and only if $Y$ is a global
 intersection.}

 \medskip

 {\bf Proof.}
The Cech cocycle $h_{ij}$ is the boundary of the $M_X^l$
$0$-cocycle $F_j-F_j$, this implies that $(h_{ij})_{i,j\in I}$ is
a $1$-Cech $Z_X^l$-cocycle. Suppose that the class $[h_{ij}]$ of
$(h_{ij})$ vanishes, this implies the existence of a $0$-chain
$f_i$ of $Z_X(l)$, such that $h_{ij}=f_j^i-f_i^j$. The boundary of
$(F_i-f_i)$ is zero, this implies that $F_i-f_i$ is the
restriction of a global section $F$ of $M_X^l$; we can suppose
that $F$ is the class of $(h_1,...h_c)$ since $Y$ is irreducible.
This implies that $Y$ is the locus of $h_1,...,h_c$.

Conversely, suppose that $Y$ is the complete intersection of
$(f_1,...,f_l)$. Then we can take $F_i$ to be the restriction to
$U_i$ of the projection of $(f_1,...,f_l)$ to $M_X$. This implies
the result.

\bigskip

{\bf Examples.}

\medskip

Suppose that $X=Spec(k)$, $L_X=k^*$ the set of non zero elements
of $k$, for every element $a\in k^*$, $(a)=0$. This implies that
$Z_X(l)={k^*}^{\otimes l}$.

\medskip

Suppose that $X$ is a curve; if $l>1$, $CH_X^l=0$.

\medskip

{\bf Proposition 2.} {\it Let $X=P^2k$ there exists a non trivial
$Z_X(2)$-gerbe defined on $X$.}

\medskip

{\bf Proof.} First we construct a non trivial element of
$H^1(X,CH_X^2)$. We can cover $X$ with the three open subsets
$U_i=\{[X_1,X_2,X_3], X_i\neq 0\}, i=1,2,3$. On $U_i\cap U_j$,
$c_{ij}$ is the homogeneous element whose $i$ and $j$ coordinates
are $1$, and the other is $0$; it is the intersection of the lines
defined by $X_i-X_j$ and $X_k, k\neq i,j$. Since
$U_{ijk}=\{[X_1,X_2,X_3]: X_i\neq 0, i=1,2,3\}$, it implies that
$c_{ij}^k=0$. Thus the family $(c_{ij})$ defines a cocycle. This
cocycle is not a boundary: Suppose that there exists a chain
$(c_i)_{i\in I}$ such that $c_{ij}=c_j^i-c_i^j$. Write
$c_i=l_i^1h_i^1+..+l_i^{j_1}h_i^{j_i},i=1,2,3,
l_i^1,..,l_i^{j_i}\in{\Z}, h_i^{j_n}\in U_i$. Suppose that the
second homogeneous coordinate of a component $h_1^{j_l}$ is not
zero, then its third homogeneous coordinate is not zero since
$c_{13}=c_3^1-c_1^3$, if its third homogeneous coordinate is not
zero, then its second homogeneous coordinate is not zero also.
Since its first coordinate is not zero, we deduce that the
coordinates of $h_1^{j_l}$ are not zero, this argument implies
that $h_i^{j_l}\in U_{ijk}, i=1,2,3$. This is in contradiction
with the fact that $c_{ij}=c_j^i-c_i^j$. Thus the class
$(c_{ij})_{i,j\in I}$ defines a non trivial $Z_{P^2k}(2)$-gerbe on
$P^2k$.

We have to show that the class $c$ defined by $(c_{ij})$ is not in
the image of $H^1(X,M_X^2)\rightarrow H^1(X,CH_X^2)$. Suppose it
is in that image, and let $h$ be an element in its preimage. We
denote by $h_{ij}$, the value of $h$ on $U_{ij}$. Suppose
$i=1,j=2$ we can represent it by a  couple $(f_{12},g_{12}), \in
L_X^2$ such that the intersection of the divisors of $f_{12}$ and
$g_{12}$ is $[1,1,0]$. Since on $U_{123}$, we can write $h^3_{12}$
as a combination of $h^2_{13}$ and $h^1_{23}$ by applying the
cocycle condition, this combination can be written in $U_{12}$,
but this is impossible, since the locus of the components of
$h_{13}$ and $h_{23}$ does not contain in $c_{12}$ $\bullet$

\bigskip

{\bf The \'Etale topology.}

\medskip

 Suppose that $X$ is an integral quasi-compact scheme
   equipped with the \'etale topology, we denote
respectively  by $U_X$ and $Div_X$ the sheaves of non zero
rational functions and the quotient of $U_X$ by $O_X^*$, the sheaf
of non zero regular functions  defined on the \'etale topology.
Let $Z_X(1)$ be the kernel of the map $U_X\rightarrow U_X/O_X^*$,
we have an exact sequence: $H_{et}^n(X,Z_X(1))\rightarrow
H^n_{et}(X,U_X)\rightarrow H^n_{et}(X,Div_X)\rightarrow
H^{n+1}_{et}(X,Z_X(1))$.

We can define the notion of $p$-Cartier \'etale gerbe to be the
quotient of $H^p(X,Div_X)$ by the image of the map $H^p(X,
U_X)\rightarrow H^p(X, Div_X)$. It is shown that if $p=1$,
$H^1_{et}(X,Z_X(1))=H_{Zar}(X,Z_X(1))=Pic(X)$.

If $X$ is smooth, then the sheaf of \'etale Cartier divisors can
be identified with the sheaf of Weil \'etale divisors whose
sections are summands of irreducible codimension $1$ varieties.
This implies that $H^i(X,Div_X)=\sum_{x\in X^1}H^i(k(x),{\Z}),
i=1,2$, where $X^1$ is a closed point of codimension $1$, and
$k(x)$ its residue field. The Hilbert 90 theorem implies that
$H^1(k(x),{\Z})=0$. This implies that the $1$-\'etale Cartier
gerbes are trivial if $X$ is smooth.

\medskip

{\bf The Brauer group.}

Let $k$ be a field, the Brauer group of $k$ is
$H^2_{et}(Spec(k),\bar k^*)$, where $\bar k$ is the algebraic
closure of $k$. Let $Gl(n,\bar k)$ be the linear group of
invertible $n$-matrices, and $PGl(n,\bar k)$ the corresponding
projective group. The exact sequence $1\rightarrow \bar
k^*\rightarrow Gl(n,\bar k)\rightarrow PGl(n,\bar k)\rightarrow
1$, induces an exact sequence $H_{et}^n(Spec(k),\bar
k^*)\rightarrow H_{et}^p(Spec(k),Gl(n,\bar k))\rightarrow
H_{et}^p(Spec(k),PGl(n,\bar k))\rightarrow
H_{et}^{p+1}(Spec(k),\bar k^*)$. If $p=1$, it is shown in Serre
[14] (proposition 9. p. 166) that every class in
$H_{et}^2(Spec(k),\bar k^*)$ are in the image of a morphism
$H_{et}^1(Spec(k),PGl(n,\bar k))\rightarrow H_{et}^2(Spec(k),\bar
k^*)$, for a given $n$, thus every element of the Brauer group is
the classifying cocycle of  a gerbe, which is the geometric
obtruction to lift a torsor on the \'etale topos $Spec(k)$.

\bigskip

{\bf References.}

\medskip

[1] Barr, M. Topos, triples and theories Springer (1985).
\smallskip

[2] Bloch, S. Algebraic cycles and higher K-theory. Advances in
Mathematics. 61 3 (1986) 267-304.
\smallskip

[3]  Brylinski, J-L. McLaughlin, D.A. The geometry of degree-$4$
characteristic classes and of line bundles on loop spaces II. Duke
Math. Journal 83 (1996) 105-139.
\smallskip

[4] E. Friedlander, A. Suslin, V. Voedvodsky. Cycles transfer and
motivic homology theories. Princeton University Press. (2000).
\smallskip

[5] Giraud, J. Cohomologie  non Ab\'elienne, Springer (1971).
\smallskip

[6] Giraud, J. M\'ethode de la descente. Bulletin de la
Soci\'et\'e Math\'ematiques de France (1964).

\smallskip
[7] Grothendieck, A. Th\'eorie des topos et cohomologie \'etale
des sch\'emas. I. Lectures Notes in Mathematics (269) 1972.
\smallskip

[8] Grothendieck, A. Rev\^etements \'etales et groupe fondamental,
Lectures Notes in Math. 224 (1971).
\smallskip

[9] Hartshorne, R. Algebraic geometry, Graduate Text in Math.
Springer (1977).
\smallskip

[10] Jackson, A. As If Summoned from the Void. Notice A.M.S
(2004).
\smallskip

[11] Jackson, Michael, A. A quotient of the set $[BG,BU(n)]$ for a
finite group $G$ of small rank. J. Pure and Applied Algebra
161-174 (2004).
\smallskip

[12] Knutson, D. Algebraic space. Lecture Notes in Math. 203
(1971).
\smallskip

[13] Laumon, G. Moret-Bailly, L. Champs alg\'ebriques. (1999)
\smallskip

[14] Serre, J-P. Corps locaux. Publications de l'institut
universitaire de Nangano. Hermann Paris (1968).
\smallskip

[15] Serre, J-P. Alg\`ebre locale et multiplicit\'és. Lectures
Notes in Math. 11 (1965)
\smallskip

[16] Toen, B. Champs alg\'ebriques. Cours.
\smallskip

[17] Tsemo, A. Non Abelian cohomology: the point of view of gerbed
towers. African Diaspora Journal of Math. 4. 67-85 (2005).
\smallskip

[18] Tsemo, A. Gerbes, $2$-gerbes and symplectic fibrations. Rocky
Journal of Math. 38. 727-777 (2008).
\smallskip

[19] Tsemo, A. Differentiable categories, gerbes and
$G$-structures, International Journal of Contemporary Mathematical
Sciences. 4. 1547-1590 (2009).

\end{document}